\documentclass[conference]{ieeeconf}
\IEEEoverridecommandlockouts
\usepackage{graphicx}
\usepackage{float}
\usepackage{array}
\usepackage{amsmath,amssymb}
\usepackage{eepic}
\usepackage{epsfig}
\usepackage{subfigure}
\usepackage{calc}
\usepackage{amssymb}
\usepackage{amstext}
\usepackage{amsmath}
\usepackage{multicol}
\usepackage{pslatex}
\usepackage{float}
\usepackage{listings}
\usepackage{cite}
\newtheorem{example}{Example}

\title{Limitations of the classical phase-locked loop analysis
\thanks{
    Kuznetsov N.V., Kuznetsova O.A., Leonov G.A.,  Neittaanmaki P., Yuldashev M.V., Yuldashev R.V.,
    Limitations of the classical phase-locked loop analysis, 
    Proceedings of International Symposium on Circuits and Systems (ISCAS), IEEE, 2015, pp. 533-536
}
}

\author{
Kuznetsov N.~V., 
Kuznetsova O.~A.,
Leonov G.~A.,
Neittaanm\"{a}ki P.,
Yuldashev M.~V.,
Yuldashev R.~V.
\thanks{
Faculty of Mathematics and Mechanics, Saint-Petersburg State University, Russia;
Dept. of Mathematical Information Technology, University of Jyv\"{a}skyl\"{a}, Finland; email: nkuznetsov239@gmail.com
}
}

\begin{document}

\maketitle
\thispagestyle{empty}
\pagestyle{empty}

\begin{abstract}                
 Nonlinear analysis of the classical phase-locked
loop (PLL) is a challenging task. In classical engineering
literature simplified mathematical models and
simulation are widely used for its study. In this work the
limitations of classical engineering phase-locked loop
analysis are demonstrated, e.g., hidden oscillations,
which can not be found by simulation, are discussed.
It is shown that the use of 
simplified dynamical models and the application of simulation may lead to
wrong conclusions concerning the operability of PLL-based
circuits.
\end{abstract}


\section{Introduction}

The Phase locked-loop (PLL) circuits were invented in the first half of the twentieth century
and nowadays are widely used in modern telecommunications and computers.
PLL is essentially a nonlinear control system and its {\it real model}
is described by a nonlinear nonautonomous system of differential equations
(\emph{mathematical model in the signal space}).
In practice, simulation and simplified mathematical models are widely used
for the analysis of PLL-based circuits \cite{Best-2007,PedersonM-2008-book,Tranter-2010-book}.

In the following it will be shown that
1) the use of simplified mathematical models
and
2) the application of non rigorous methods of analysis (e.g., simulation)
may lead to wrong conclusions concerning
the operability of \emph{real model} of classical PLL.

\section{Simulation of the classical phase-locked loop in MatLab Simulink}

Consider the classical PLL nonlinear models in the signal
and signal's phase spaces \cite{Viterbi-1966,Gardner-1966,LeonovKYY-2012-TCASII,LeonovKYY-2015-SIGPRO,LeonovK-2014,BestKLYY-2014-IJAC}. \begin{itemize}
  \item
    \emph{Real model of the classical PLL in the signal space} (Fig.~\ref{pll-signal-ics})
    or its \emph{nonlinear mathematical model in the signal space}
    (corresponds to the SPICE-level simulation):

    \begin{equation}\label{sysPLL-signal}
     \begin{aligned}
     & \dot{x} = A x + b\varphi(t), \quad \varphi(t) =  \sin(\theta_1(t))\cos(\theta_2(t)) \\
     & \dot\theta_{1} \equiv \omega_1, \ \dot\theta_{2} = \omega_2^{\text{free}} + L(c^*x) + Lh\varphi(t).\\
     \end{aligned}
    \end{equation}

    \begin{figure}[thpb]
    \centering
    \includegraphics[width=0.45\textwidth]{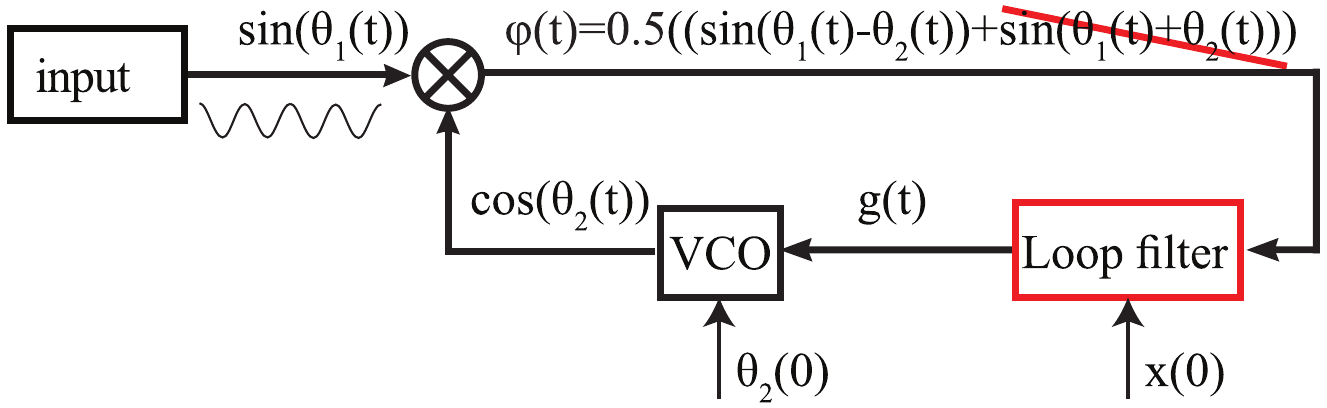}
    \caption{Real model of the classical PLL in the signal space}
    \label{pll-signal-ics}
    \end{figure}

  \item
    \emph{Model of the classical PLL in signal's phase space}
    (Fig.~\ref{pll-phase-ics});
    system \eqref{sysPLL-signal} with averaged $\varphi(t)\approx\varphi(\theta_{\Delta}(t))$
    gives the \emph{nonlinear mathematical model in signal's phase space}:
    \begin{equation}\label{sysPLL-phase}
     \begin{aligned}
     & \dot{x} = A x + b \varphi(\theta_{\Delta}), \\
     & \dot\theta_{\Delta} = \omega_{\Delta} - L(c^*x) - Lh\varphi(\theta_{\Delta}),\\
     & \theta_{\Delta}(t)=\theta_1(t)-\theta_2(t), \ \omega_{\Delta} \equiv \omega_1-\omega_2^{\text{free}}.
     \end{aligned}
    \end{equation}

    \begin{figure}[thpb]
    \centering
     \includegraphics[width=0.45\textwidth]{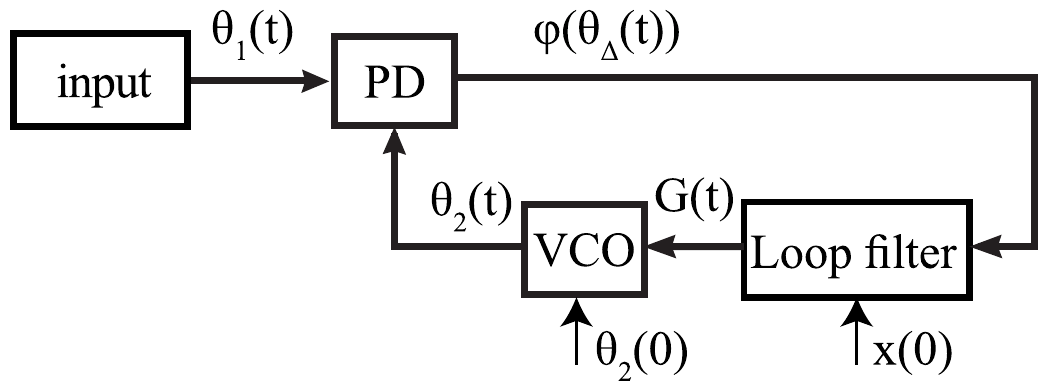}
    \caption{Simplified model of the classical PLL in signal's phase space}
    \label{pll-phase-ics}
    \end{figure}
\end{itemize}

\begin{figure}[H]
 \includegraphics[width=0.9\linewidth]{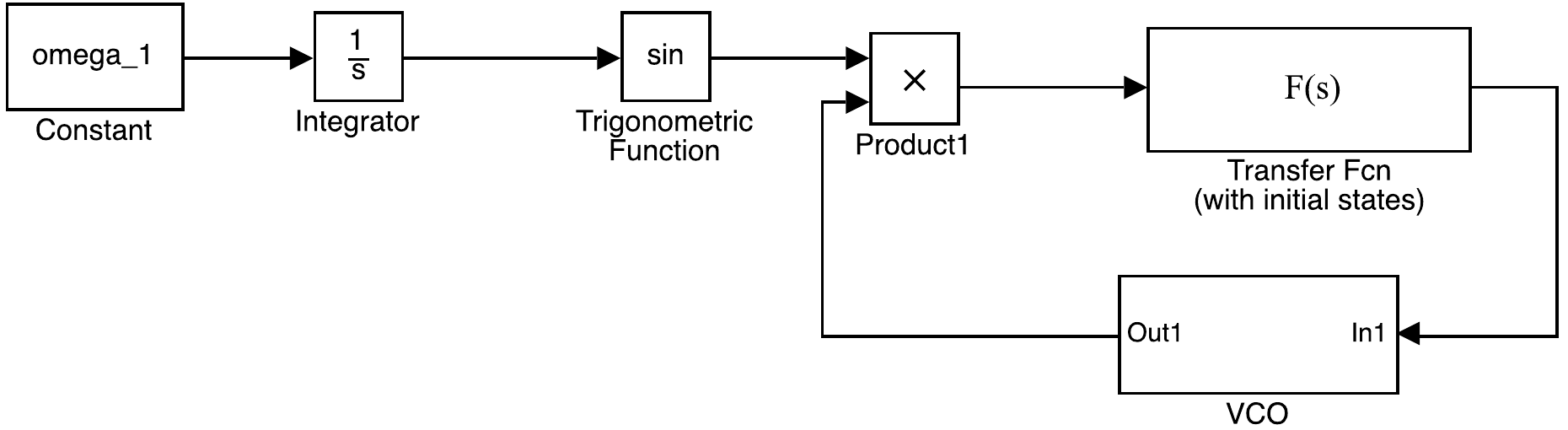}
 \centering
\caption{Simulink realization of the real model in the signal space}
\label{simulink_model_signal}
\end{figure}

\begin{figure}[H]
 \includegraphics[width=0.9\linewidth]{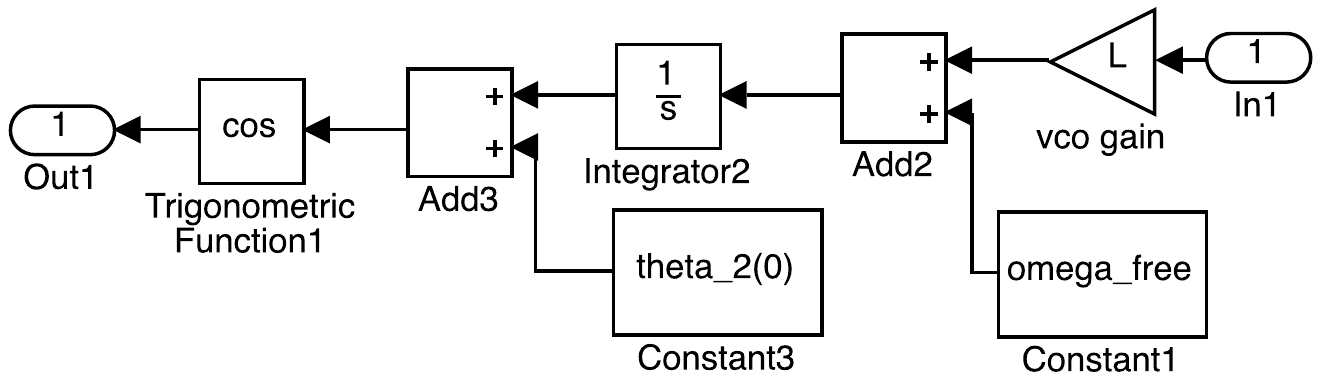}
 \centering
\caption{Simulink realization of the VCO for the real model}
\label{simulink_vco_signal}
\end{figure}

Let us construct MatLab Simulink model, which corresponds to
the  model in the signal space (see Fig.~\ref{simulink_model_signal}).
\vspace{-0.2in}
\noindent Here all elements are standard blocks from Simulink Library
except for the VCO.
The VCO subsystem is shown in Fig.~\ref{simulink_vco_signal}.
\vspace{-0.2in}
\noindent The VCO subsystem consists of one input, which is amplified by
L (Gain block).
 The integration of the sum of amplified input signal
 and the VCO free-running frequency omega\_free
forms the phase of the VCO output. The VCO output corresponds to $\cos(\cdot)$.

\begin{figure}[H]
 \includegraphics[width=0.85\linewidth]{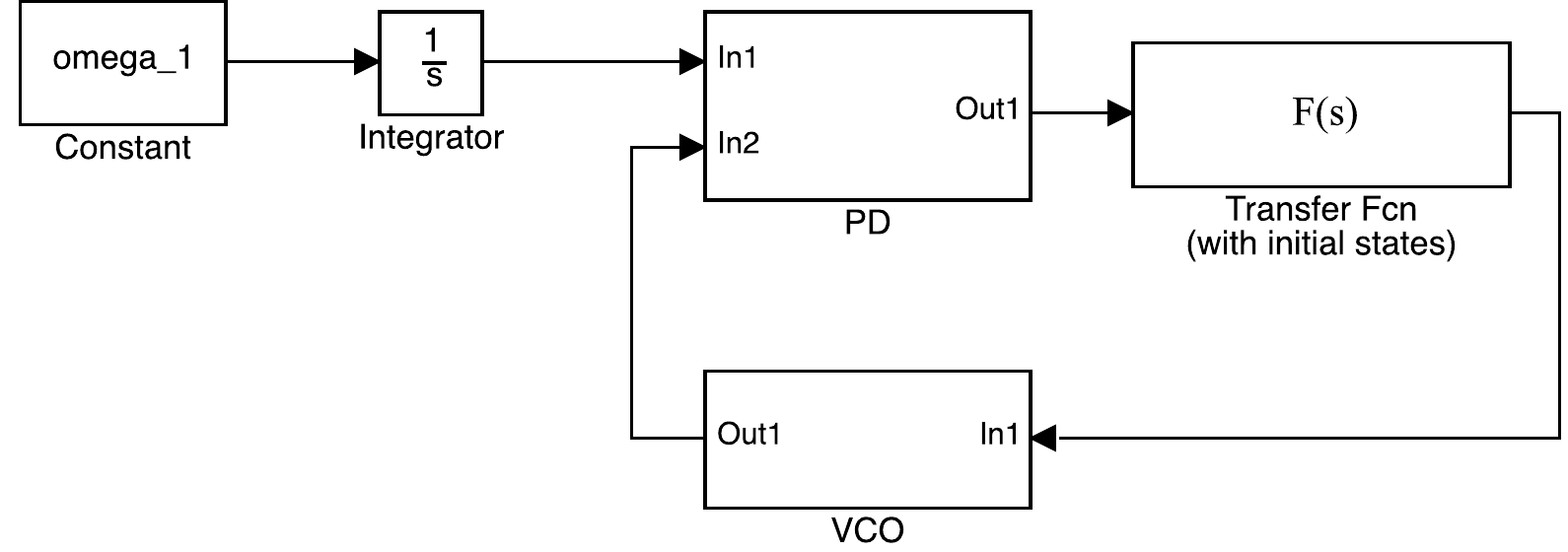}
 \centering
\caption{Simulink realization of the model in signal's phase space}
\label{simulink_model_phase}
\end{figure}

Now consider MatLab Simulink model, which corresponds to
 the model in signal's phase space (see Fig.~\ref{simulink_model_phase}).

\begin{figure}[H]
 \includegraphics[width=0.85\linewidth]{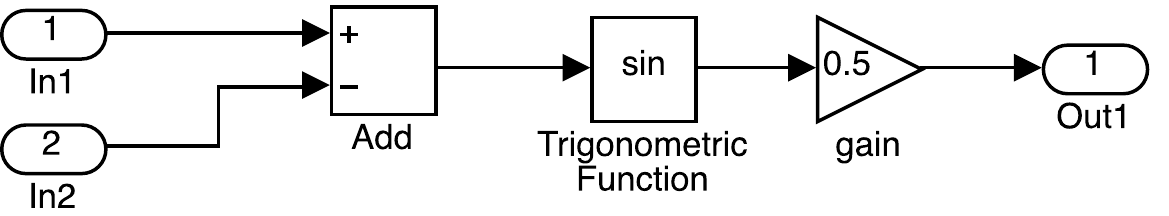}
 \centering
\caption{Simulink realization of the PD in signal's phase space}
\label{simulink_PD_phase}
\end{figure}

The PD subsystem is shown in Fig.~\ref{simulink_PD_phase}.

\begin{figure}[H]
 \includegraphics[width=0.85\linewidth]{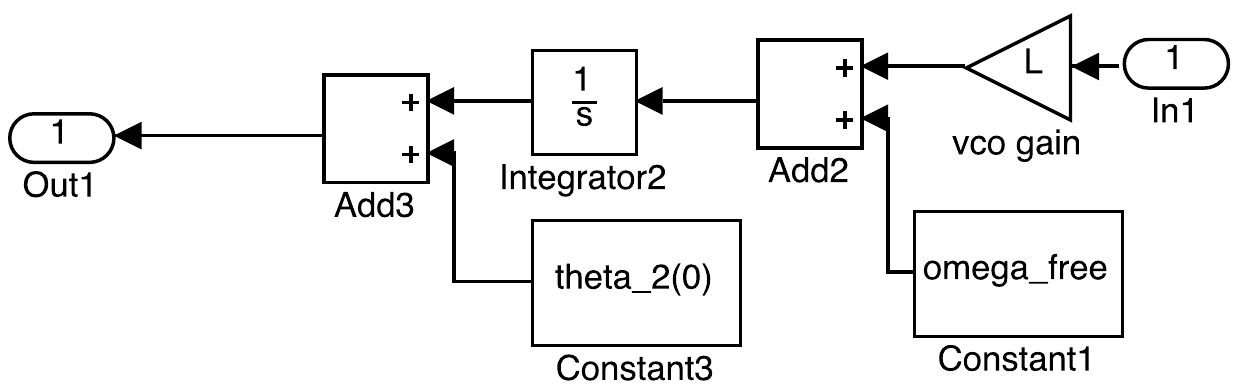}
 \centering
\caption{Simulink realization of the VCO in signal's phase space}
\label{simulink_vco_phase}
\end{figure}

\noindent The VCO subsystem in signal's phase space is shown in Fig.~\ref{simulink_vco_phase}.
The VCO output in signal's phase space corresponds to $\theta_2(t)$.

\subsection{Simulation parameters and examples}
Consider a passive lead-lag loop filter
with the transfer function
$F(s) = \frac{1+s \tau_2}{1+s(\tau_1 + \tau_2)}$, $\tau_1 = 0.0448$, $\tau_2 = 0.0185$
and the corresponding parameters
$A = -\frac{1}{\tau_1+\tau_2}$,
$b = 1 - \frac{\tau_2}{\tau_1+\tau_2}$,
$c = \frac{1}{\tau_1+\tau_2}$,
$h = \frac{\tau_2}{\tau_1+\tau_2}$.
The input signal frequency is $\omega_1=100000$,
initial phase is zero: $\theta_1(0)=0$,
and the VCO input gain $L = 250$.

\begin{example}\rm
This example shows the importance of initial state of filter (see Fig.~\ref{lpf_ics}):
while the real model (see Fig.~\ref{pll-signal-ics})
with nonzero initial state of loop filter $x_0 = 0.18$ does not acquire lock (black color),
the same real model with zero initial state of loop filter $x(0) = 0$
acquires lock (red color). Here the VCO free-running frequency
$\omega_2^{\text{free}}= 100000-95$.
\begin{figure}
  \includegraphics[width=0.8\linewidth]{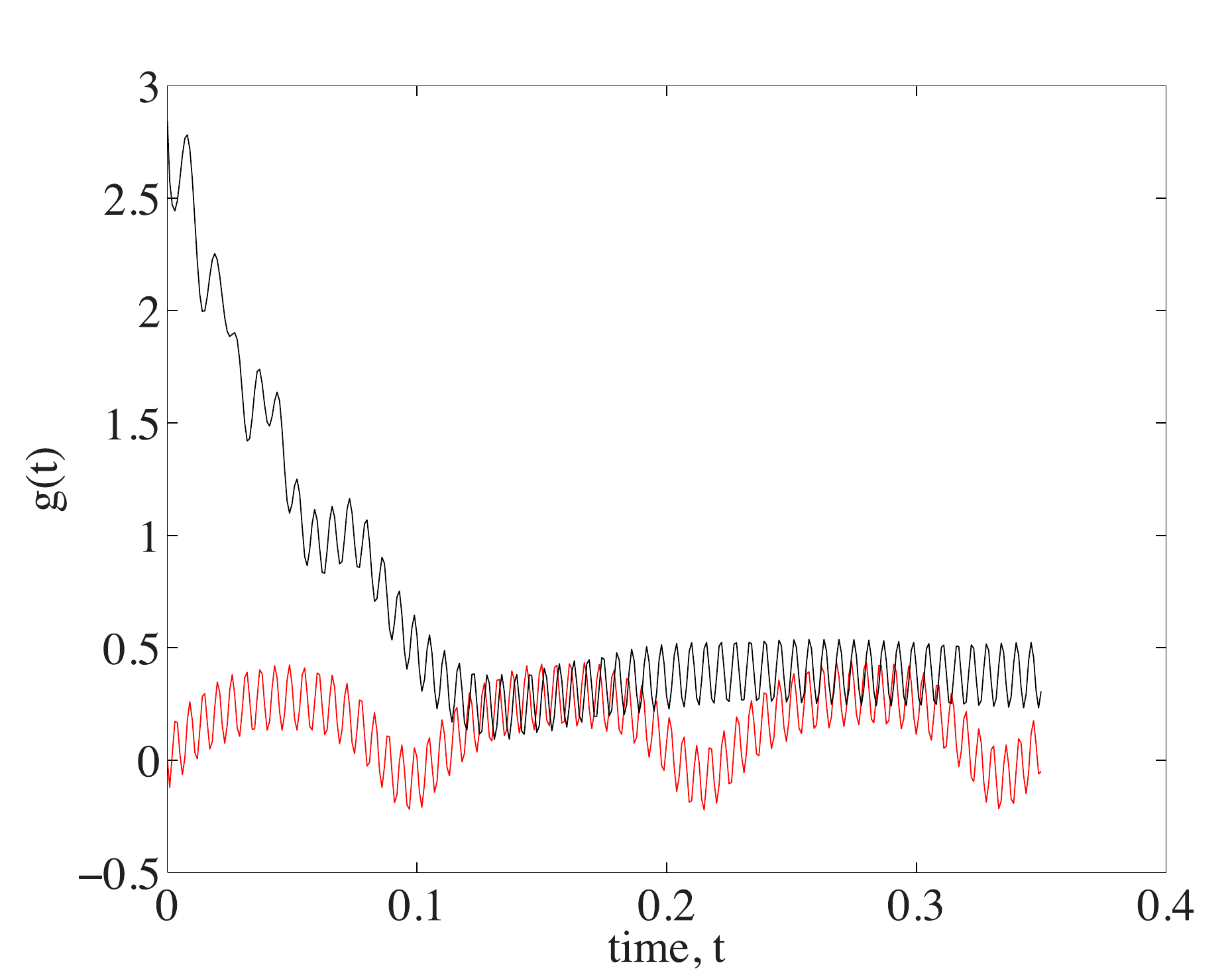}
  \caption{
  Loop filter output $g(t)$ for
  real model with nonzero initial state of loop filter (red),
  real model with zero initial state of loop filter (black).}
  \label{lpf_ics}
\end{figure}
\end{example}

\begin{example}\rm
This example shows
that the initial phase difference $\theta_1(0)-\theta_2(0)$
between the VCO signal and input signal
may affect stability of the classical PLL.
In Fig.~\ref{ex_phase} the real model (see Fig.~\ref{pll-signal-ics})
with zero initial phase difference acquire lock (red color),
the same real model with nonzero initial phase difference $\theta_{\Delta}(0) = \pi$
is out of lock (black color).
Here the VCO free-running frequency $\omega_2^{\text{free}}= 100000-95$
and the initial state of loop filter is $x_0 = 0.01$.
\begin{figure}
  \includegraphics[width=0.8\linewidth]{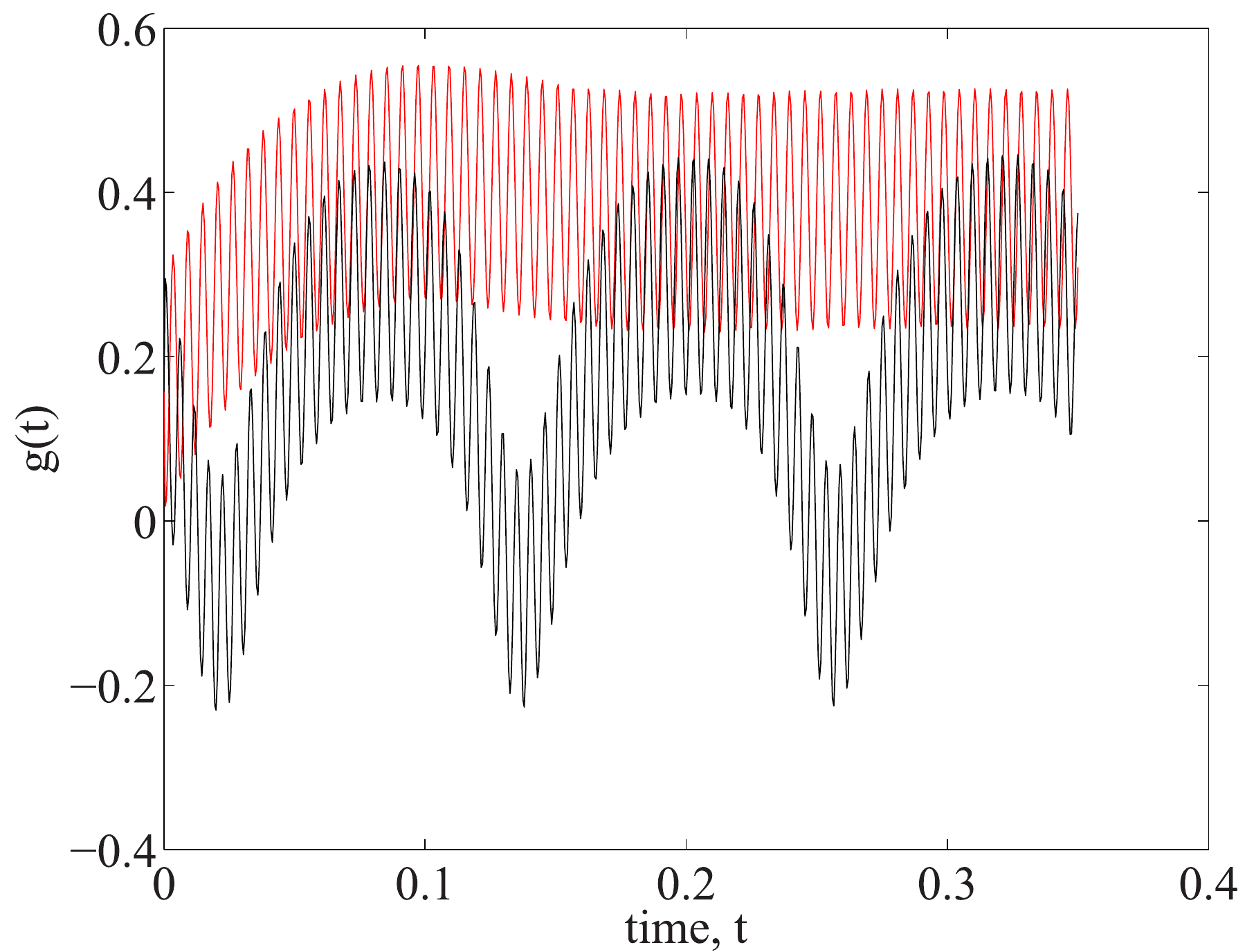}
  \caption{
  Loop filter output $g(t)$ for
  real model with nonzero initial phase difference (black),
  real model with zero initial phase difference (red).}
  \label{ex_phase}
\end{figure}
\end{example}

 Examples 1 and 2 shows that
 while the term ``initial frequency''  (without an explanation)
 is sometimes used instead of the the term ``free-running frequency''
 in engineering definitions of various stability ranges, it may lead to a misunderstanding 
 (see corresponding discussion in \cite{KuznetsovLYY-2015-IFAC-Ranges,KuznetsovLYY-2015-arxiv-Ranges}).

\begin{example}\rm
This example shows that the PLL model in signal's phase space
may not be equivalent to the PLL real model in the signal space.
In Fig.~\ref{double_freq} the real model (see Fig.~\ref{pll-signal-ics})
does not acquire lock (red color),
 the equivalent signal's phase space model acquires lock (black color).
Here the VCO free-running frequency $\omega_2^{\text{free}}= 100000-95$,
the initial state of loop filter is $x_0 = 0.017$,
 and the initial phase difference $\theta_{\Delta}(0)=2.276$.
\begin{figure}
  \includegraphics[width=0.8\linewidth]{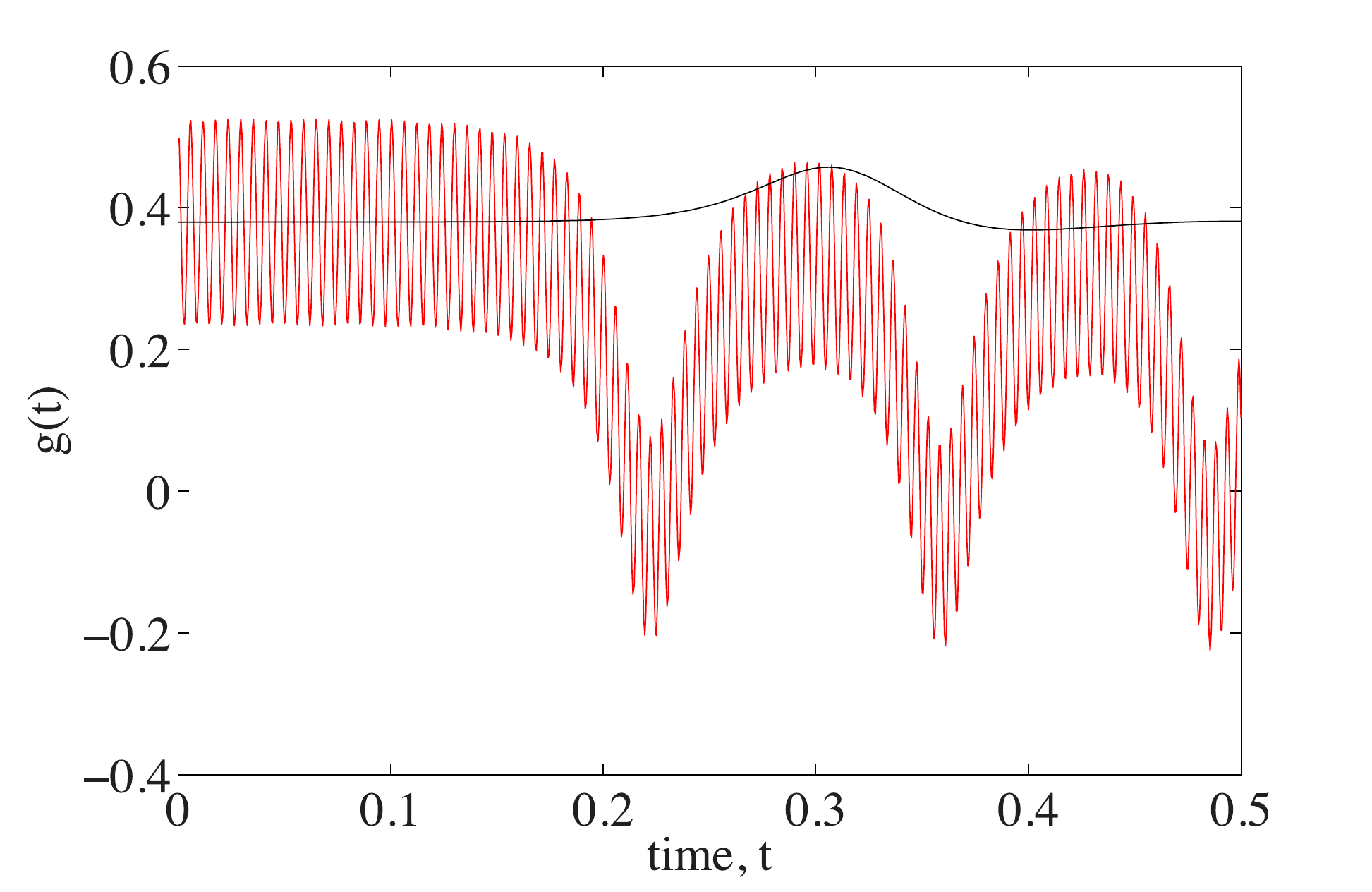}
  \caption{
  Loop filter output $g(t)$ for
  signal's phase model (black),
  real model (red).}
  \label{double_freq}
\end{figure}
\end{example}

\begin{example}\rm
 These examples shows the importance of analytic methods
for investigation of PLL stability.
More precisely, it is shown that the simulation may lead to wrong results.
In Fig.~\ref{pll_hidden}
the PLL model in signal's phase space simulated with relative tolerance ``1e-3''
does not acquire lock (black color),
but the PLL model in signal's phase space simulated with standard parameters
(relative tolerance set to ``auto'') acquires lock (red color)\footnote{
See, e.g., the corresponding internal time step parameter in PSpice [http://www.stuffle.net/references/PSpice\_help/tran.html].
In \cite{BianchiKLYY-2015}
the SIMetrics SPICE model of the two-phase PLL
with lead-lag filter 
gives two essentially different results with default sampling step and minimum sampling step set to $1m$.
}.
Here the input signal frequency is $10000$, the VCO free-running frequency $\omega_2^{\text{free}}= 10000 - 178.9$,
the VCO input gain is $L=500$,
 the initial state of loop filter is $x_0 = 0.1318$,
and the initial phase difference is $\theta_{\Delta}(0) = 0$.
\begin{figure}
  \includegraphics[width=0.85\linewidth]{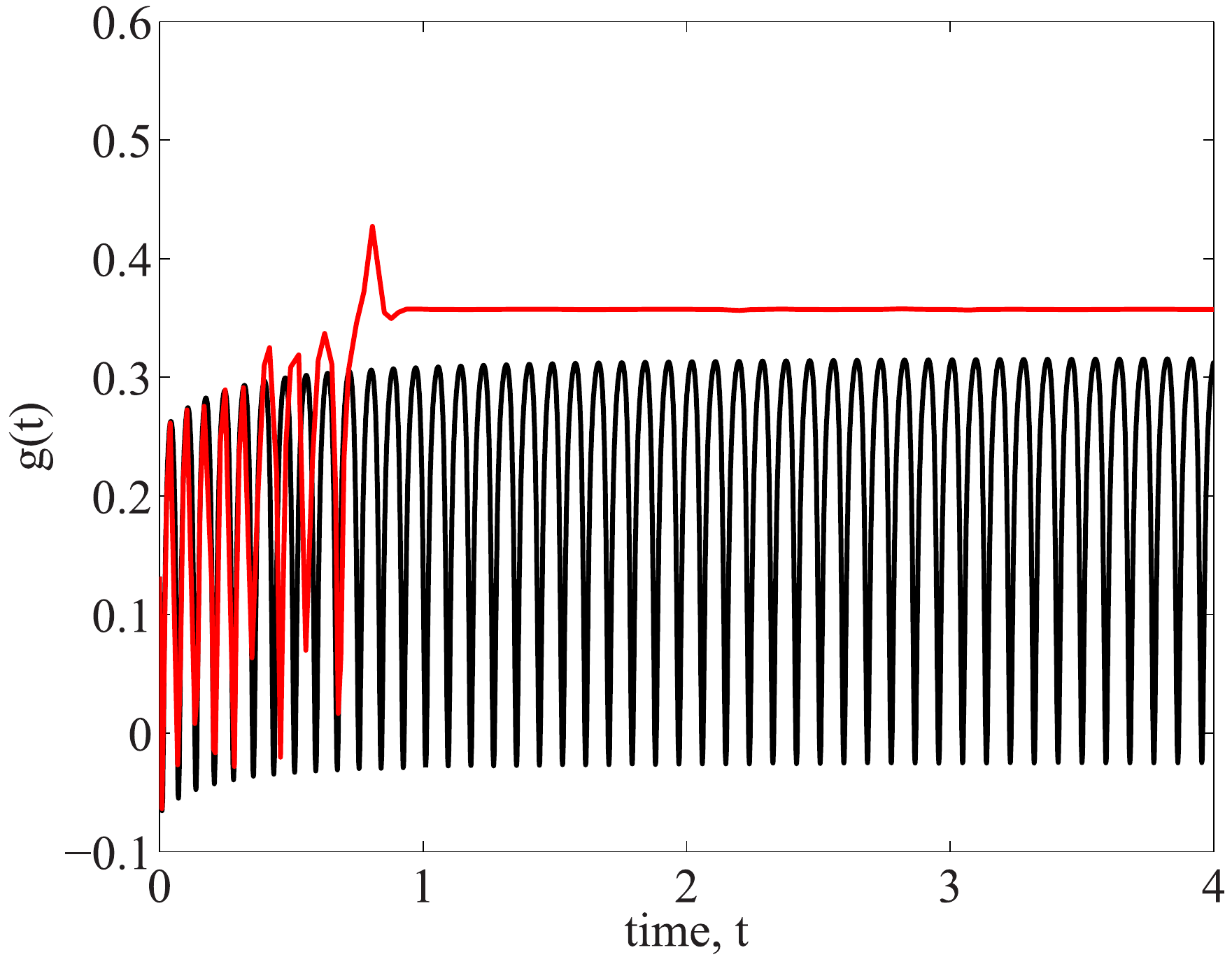}
  \caption{
  Loop filter output $g(t)$ for
  signal's phase space model with standard integration parameters (red),
  signal's phase space model with relative tolerance set to ``1e-3''(black).}
  \label{pll_hidden}
\end{figure}
Consider now a phase portrait (the loop filter state $x$
versus the phase difference $\theta_{\Delta}$)
corresponding to signal's phase model (see Fig.~\ref{pll_hidden_phase_portret}).
\begin{figure}
  \includegraphics[width=0.85\linewidth]{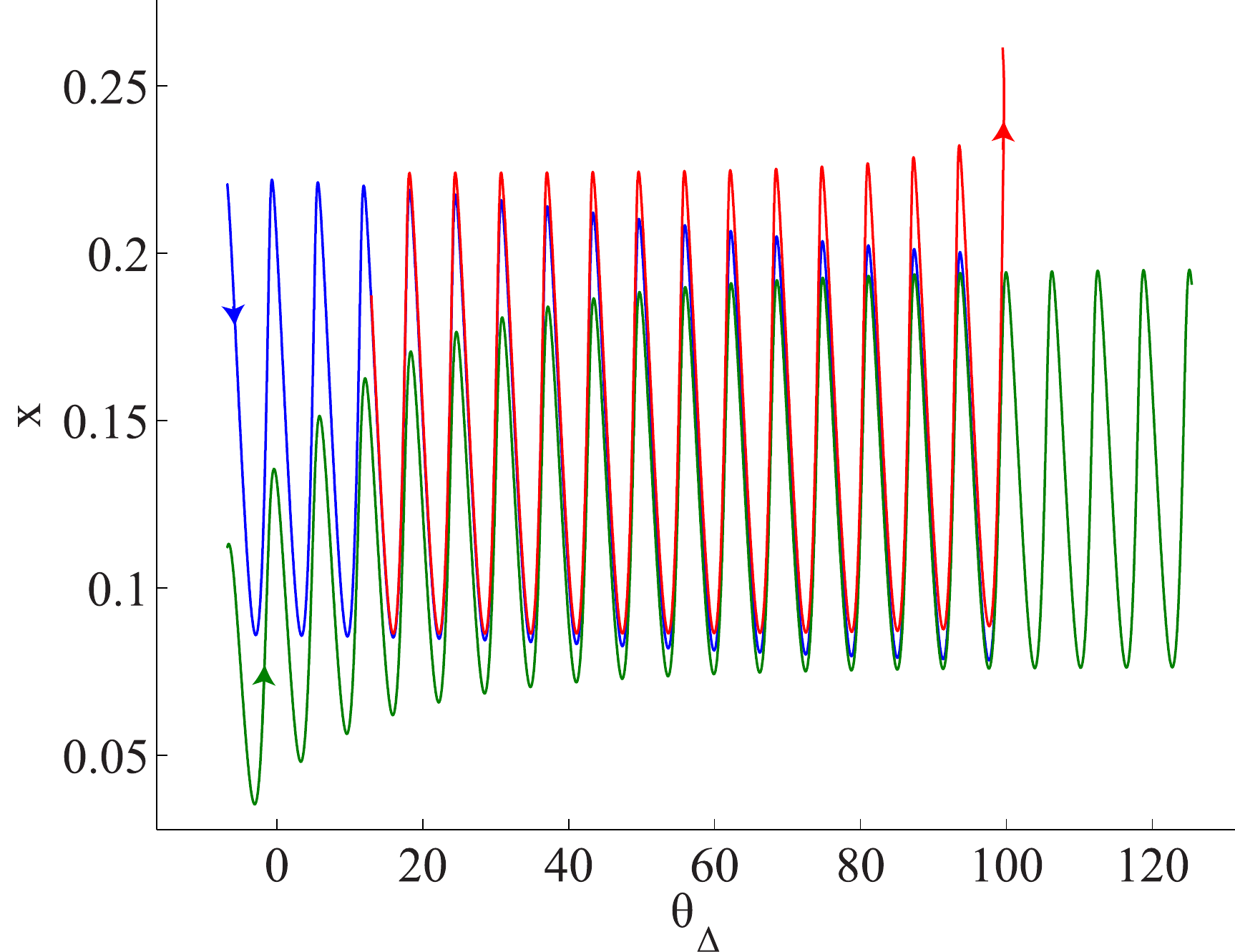}
  \caption{Phase portrait of the classical PLL with stable and unstable periodic trajectories}
  \label{pll_hidden_phase_portret}
\end{figure}
The solid blue line in Fig.~\ref{pll_hidden_phase_portret} corresponds
to the trajectory with the loop filter initial state $x(0) = 0.2206$
and the VCO phase shift $-6.808$ rad.
This line tends to the periodic trajectory,
therefore it will not acquire lock.

The solid red line corresponds to the trajectory
with the loop filter initial state $x(0) = 0.187386698333130$
and the VCO initial phase $12.938118990628919$.
This trajectory lies just under the unstable periodic trajectory and tends
to  a stable equilibrium.
In this case PLL acquires lock.

All trajectories between stable and unstable periodic trajectories
tend to the stable one (see, e.g., a solid green line).
Therefore, if the gap between stable and unstable periodic trajectories
is smaller than the discretization step,
the numerical procedure may slip through the stable trajectory.
In other words, the simulation will show that the PLL acquires lock,
but in reality it is not the case.
The considered case corresponds
to the coexisting attractors (one of which is so-called hidden oscillation)
and the bifurcation of birth of semistable trajectory \cite{LeonovK-2013-IJBC,KuznetsovLYY-2014-IFAC}.
\end{example}


An oscillation in a dynamical system can be easily localized numerically
if the initial conditions from its open neighborhood lead to long-time behavior
that approaches the oscillation.
Thus, from a computational point of view, it is natural to suggest
the following classification of attractors,
based on the simplicity of finding the basin of attraction in the phase space \cite{KuznetsovLV-2010-IFAC,LeonovKV-2011-PLA,LeonovKV-2012-PhysD,LeonovK-2013-IJBC}:
{\it An attractor is called a \emph{hidden attractor} if its
 basin of attraction does not intersect with
 small neighborhoods of equilibria,
 otherwise it is called a \emph{self-excited attractor}.
}

For a \emph{self-excited attractor} its basin of attraction
is connected with an unstable equilibrium
and, therefore, self-excited attractors
can be localized numerically by the
\emph{standard computational procedure},
in which after a transient process a trajectory,
started from a point of an unstable manifold in
a neighborhood of an unstable equilibrium,
is attracted to the state of oscillation and traces it.
Thus self-excited attractors can be easily visualized.

In contrast, for a hidden attractor its basin of attraction
is not connected with unstable equilibria.
For example, hidden attractors are attractors in
the systems with no equilibria
or with only one stable equilibrium
(a special case of multistable systems and
coexistence of attractors).

\section{Conclusion}

The derivation of
mathematical model in signal's phase space
and the use of the results of its analysis
to draw conclusions about the behavior
of real model in the signal space
 have need for a rigorous foundation.
But the attempts to justify analytically the reliability of conclusions,
based on such engineering approaches,
and to study the nonlinear models of PLL-based circuits are quite rare
in the modern engendering literature \cite{Abramovitch-2002}.
One of the reasons is that
``{\it nonlinear analysis techniques are well beyond the scope of
most undergraduate courses in communication theory}'' \cite{Tranter-2010-book}.

 The examples considered in the paper are the motivation
to use rigorous analytical methods
for the analysis of nonlinear PLL models \eqref{sysPLL-signal}-\eqref{sysPLL-phase}.
Some analytical tools can be found in \cite{GeligLY-1978,LeonovRS-1992,Stensby-1997,SuarezQ-2003,Margaris-2004,KudrewiczW-2007,LeonovK-2014}.

Note once more that various simplifications
and the analysis of linearized models of control systems
may result in incorrect conclusions (see, e.g.,
  the counterexamples to the filter hypothesis,
Aizerman's and Kalman's conjectures on the absolute stability of nonlinear
control systems \cite{BraginVKL-2011,LeonovK-2013-IJBC},
and the Perron effects of the largest Lyapunov exponent sign reversals \cite{KuznetsovL-2005}, etc.).

In the work it is shown that
1) the consideration of simplified models,
constructed intuitively by engineers and
2) the application of non-rigorous methods of analysis
(e.g., simulation and linearization)
can lead to wrong conclusions concerning the operability
of the classical phase-locked loop.
Similar examples for nonlinear Costas loop models can be found in \cite{KuznetsovKLNYY-2014-ICUMT-QPSK,KuznetsovKLSYY-2014-ICUMT-BPSK,KuznetsovL-2014-IFACWC,KudryashovaKKLSYY-2014-ICINCO,BestKKLYY-2015-ACC}.

\section*{\uppercase{Acknowledgements}}
 Authors were supported by Saint-Petersburg State University
 and Russian Scientific Foundation.
 The authors would like to thank Roland~E.~Best
 (the founder of Best Engineering company, Switzerland;
 the author of the bestseller on PLL-based circuits \cite{Best-2007})
 for valuable discussions.

\bibliographystyle{IEEEtran}
\newcommand{\noopsort}[1]{} \newcommand{\printfirst}[2]{#1}
  \newcommand{\singleletter}[1]{#1} \newcommand{\switchargs}[2]{#2#1}

\end{document}